\documentclass[titlepage,12pt]{article}

\usepackage{graphicx}
\usepackage{amssymb}

\newcommand{\R}{\mathbb R}
\newcommand{\C}{\mathbb C}

\newtheorem{teo}{Theorem}[section]
\newtheorem{lema}[teo]{Lemma}

\newtheorem{cor}[teo]{Corollary}
\newtheorem{remark}[teo]{Remark}

\begin{document}

\baselineskip=18pt

\begin{center}
{\Large{\bf Hopf Bifurcation in a Model for \\ Biological Control }}
\end{center}

\begin{center}
{\large Jorge Sotomayor \\} {\em Instituto de Matem\'atica e
Estat\'{\i}stica, Universidade de S\~ao Paulo\\ Rua do Mat\~ao
1010, Cidade Universit\'aria, \\ CEP 05.508-090, S\~ao Paulo, SP,
Brazil\\}e--mail: sotp@ime.usp.br
\end{center}

\begin{center}
{\large Luis Fernando Mello \\}
{\em Instituto de Ci\^encias
Exatas, Universidade Federal de Itajub\'a\\Avenida BPS 1303,
Pinheirinho, CEP 37.500-903, Itajub\'a, MG, Brazil
\\}e--mail: lfmelo@unifei.edu.br
\end{center}

\begin{center}
{\large Danilo Braun Santos \\} {\em Centro de Ci\^encias Sociais
e Aplicadas, Universidade Mackenzie \\ Rua Itamb\'e, 45, Consola\c
c\~ao, CEP 01302-907, S\~ao Paulo, SP, Brazil
\\}e--mail: danilobraun@mackenzie.com.br
\end{center}

\begin{center}
{\large Denis de Carvalho Braga \\}
{\em Instituto de Sistemas
El\'etricos e Energia, Universidade Federal de Itajub\'a\\Avenida
BPS 1303, Pinheirinho, CEP 37.500-903, Itajub\'a, MG, Brazil
\\}e--mail: braga@unifei.edu.br
\end{center}

\begin{center}
{\bf Abstract}
\end{center}

\vspace{0.1cm}

In this paper we study the Lyapunov stability and Hopf bifurcation
in a biological system which models the biological control of
parasites of orange plantations.

\vspace{0.1cm}



\newpage
\section{\bf Introduction of the Mathematical Model}\label{intro}

In this work we study a system of four coupled differential
equations (\ref{sistemafinal}) which models the interaction
between two biological species, each presenting two stages in
their metamorphosis, living in a common habitat with limited
resources.

The differential equations analyzed here are
\begin{eqnarray}\label{sistemafinal}
P'&=&\frac{dP}{dt}=\phi_1 \left( 1-\frac{M}{c_1} \right) M-(\alpha_1 + \beta_1 )P - k_1 P G \nonumber \\
M'&=&\frac{dM}{dt}=\alpha_1 P - \mu_1 M   \\
L'&=&\frac{dL}{dt}=\phi_2 \left( 1-\frac{G}{c_2} \right) G - (\alpha_2 + \beta_2)L + k_2 P G \nonumber  \\
G'&=&\frac{dG}{dt}=\alpha_2 L - \mu_2 G.\nonumber
\end{eqnarray}

This model ---an elaboration of Lotka-Volterra equations, taking
into account the stages or compartments in the biological
populations--- was proposed by Yang and Ternes \cite{hyu, hyuson}
and Ternes \cite{son} for a study of the biological
control\footnote{http://en.wikipedia.org/wiki/Biological\_control}
of orange plantations leaf parasites $P$, which is a pre-adult
stage for $M$, by their natural enemies $L$, which is an early
stage for $G$.

Other differential equations have been proposed as models for
interacting populations partitioned in compartments, representing
several situations of biological interest. See, among many others,
Hethcote et al. \cite{HH}, Jacquez and Simon \cite{cs} and Godfray
and Waage \cite{gw}.

In \cite{hyu, hyuson} and \cite{son} $P$ and $M$ are the densities
of pupae\footnote{http://en.wikipedia.org/wiki/Pupa} and female
adults of {\it Phyllocnistis citrella} (which in its
larva\footnote{http://en.wikipedia.org/wiki/Larva} stage is the
citrus leafminer\footnote{http://www.agrobyte.com.br/minadora.htm;
\, http://en.wikipedia.org/wiki/Citrus}), $L$ and $G$ are the
densities of larvae and female adults of its native parasitoid
{\it Galeopsomyia fausta} (whose larvae feed on the pupae of
$M$\footnote{http://www.seea.es/conlupa/AlbertoWeb/framesparasitoides.htm.
This site has impressive photos of hosts and parasitoids.}).

The meaning of the parameters in (\ref{sistemafinal}), where the
notation of \cite{hyu, son} has been preserved, is as follows:
$\alpha_1$ is the rate of pupae that give rise to adults $M$,
$\beta_1$ is the mortality rate of pupae, $\mu_1$ is the mortality
rate of adults $M$, $\phi_1$ is the rate of eggs that give rise to
pupae, $c_1$ is the carrying capacity of the population $M$,
$\alpha_2$ is the rate of larvae that, evolving through pupae,
give rise to adults $G$, $\beta_2$ is the mortality rate of larvae
and pupae, $\mu_2$ is the rate of mortality of adults $G$,
$\phi_2$ is the
oviposition\footnote{en.wikipedia.org/wiki/Oviposition} rate of
the parasite and $c_2$ is the carrying capacity of the population
$G$. Here we assume that the pupa (respectively larva) population
decreases (respectively  increases) at a rate proportional to
$P-G$ encounters that is $k_1 P G$ (respectively  $k_2 P G$).

This model represents the evolution of female populations. If
necessary the male populations can be estimated using the sexual
ratio of each species.

\begin{remark} All the parameters $\alpha_1 ,
\beta_1 , \mu_1 , \phi_1 , c_1, k_1, \alpha_2, \beta_2, \mu_2,
\phi_2, c_2, k_2$ are positive. As the damage to the $P$
population must be larger than the benefit to the $L$ population
it is natural to assume that $k_1 \geq k_2$. \label{parametros}
\end{remark}

Here will be established the location and the stability character
of the equilibria of (\ref{sistemafinal}), four in number. Also is
determined the bifurcation variety in the space of parameters,
representing the transition from asymptotically stable to saddle
type at the equilibrium point with positive coordinates,
representing the coexistence of the two species. See Theorem
\ref{teo2} and its Corollary \ref{hopf}.

Fixing the all the parameters in (\ref{sistemafinal}) to
biologically feasible values, taken from \cite{son} and
\cite{dan}, but letting the interaction coefficients $k_1$ and
$k_2$ vary in a positive quadrant, the nature of the bifurcation
phenomenon in this plane by crossing the bifurcation curve is
established. See Theorem \ref{teoremahopf} and Figure
\ref{admpara}. This is done by means of a computer assisted
calculation of the first Lyapunov coefficient, found to be
positive. The Hopf bifurcation analysis of this point implies that
the bifurcating periodic orbit is asymptotically unstable, of
saddle type which surrounds an attracting equilibrium with small
attracting basin. The dependence of the bifurcation curve on the
parameter $c_2$ is studied in Theorem \ref{moves} and illustrated
in Figure \ref{curvasigma}.

In Section \ref{conclusion} the implications of the results in
this paper are discussed and interpreted from a wider perspective.

\section{ Stability Analysis of Equilibria}\label{stability}

Assume the following notation:
\begin{equation} \label{R1R2}
R_1 = \frac{\alpha_1 \phi_1 }{\mu_1 (\alpha_1  + \beta_1 )}, \:\:
R_2 = \frac{\alpha_2 \phi_2} {\mu_2 (\alpha_2 + \beta_2)}.
\end{equation}

The differential equations (\ref{sistemafinal}) have four
equilibria
\begin{equation}
\mathcal A_1 = (P_1,M_1,L_1,G_1)=(0,0,0,0), \label{ponto1}
\end{equation}
\begin{equation}
\mathcal A_2 = (P_2,M_2,L_2,G_2)=\left( \frac{c_1 \mu_1}{\alpha_1}
\left( 1-\frac{1}{R_1} \right),c_1 \left( 1-\frac{1}{R_1}
\right),0,0 \right), \label{ponto2}
\end{equation}
\begin{equation}
\mathcal A_3 = (P_3,M_3,L_3,G_3)=\left( 0,0,\frac{c_2
\mu_2}{\alpha_2} \left( 1-\frac{1}{R_2} \right),c_2 \left(
1-\frac{1}{R_2}\right) \right),\label{ponto3}
\end{equation}
and
\begin{equation}
\mathcal A_4 = (P_4,M_4,L_4,G_4), \label{ponto4}
\end{equation}
where
\[
P_4 = \frac{c_1 \mu_1 \phi_2}{\alpha_1 ^2 \phi_1 \phi_2 + \mu_1 ^2
c_1 c_2 k_1 k_2} \left(\alpha_1 \phi_1 \left( 1-\frac{1}{R_1}
\right) -\mu_1 c_2 k_1 \left( 1-\frac{1}{R_2}\right) \right),
\]
\[
M_4 =  \frac{c_1 \alpha_1 \phi_2}{\alpha_1 ^2 \phi_1 \phi_2 +
\mu_1 ^2 c_1 c_2 k_1 k_2} \left( \alpha_1 \phi_1 \left(
1-\frac{1}{R_1} \right)- \mu_1 c_2 k_1 \left( 1-\frac{1}{R_2}
\right) \right),
\]
\[
L_4 = \frac{c_2 \mu_2 \alpha_1 \phi_1}{\alpha_2 (\alpha_1 ^2
\phi_1 \phi_2 + \mu_1 ^2 c_1 c_2 k_1 k_2)} \left( c_1 \mu_1 k_2
\left( 1-\frac{1}{R_1} \right) +\alpha_1 \phi_2 \left(
1-\frac{1}{R_2} \right) \right),
\]
\[
G_4 = \frac{c_2 \alpha_1 \phi_1}{\alpha_1 ^2 \phi_1 \phi_2 + \mu_1
^2 c_1 c_2 k_1 k_2} \left( c_1 \mu_1 k_2 \left(
1-\frac{1}{R_1}\right) + \alpha_1 \phi_2 \left( 1-\frac{1}{R_2}
\right) \right).
\]

\begin{remark} \label{remark1}
If $R_1 > 1$ and $R_2 > 1$, then the equilibria $\mathcal A_1$,
$\mathcal A_2$ and $\mathcal A_3$, have only non-negative
coordinates. If $k_1 < k_{1_{max}}$, where
\begin{equation}
k_{1_{max}} = \frac{\alpha_1 \phi_1 \left(
1-\frac{1}{R_1}\right)}{c_2 \mu_1 \left( 1-\frac{1}{R_2}\right)},
\label{k1max}
\end{equation}
then the coordinates of the equilibrium $\mathcal A_4$ are also
non-negative.
\end{remark}

The Jacobian matrix of (\ref{sistemafinal}) at ${\bf x}=(P,M,L,G)
\in \R ^4$ has the form
\begin{equation}\label{jacob}
J({\bf x})=\left(\begin{array}{cccc}-\alpha_1 - \beta_1 -k_1 G &
\phi_1
- \frac{2 \phi_1 M}{c_1} & 0 & -k_1 P \\
\alpha_1 & - \mu_1 & 0 & 0\\
k_2 G & 0 & - \alpha_2 - \beta_2 & \phi_2 - \frac{2 \phi_2 G}{c_2}+ k_2 P \\
0 & 0 & \alpha_2 & - \mu_2
\end{array}\right),
\end{equation}
while its characteristic polynomial is given by
\begin{equation}
p (\lambda) = det (J({\bf x}) - \lambda I) = \Theta_1 \: \Theta_2  \\
+ \alpha_2 (\mu_1 + \lambda) k_1 k_2 P G , \label{caracteristica}
\end{equation}
where
\[
\Theta_1 = (\mu_1 + \lambda)(\alpha_1 +\beta_1 +k_1 G+
\lambda)-\alpha_1 \left( \phi_1 - \frac{2 \phi_1 M}{c_1} \right)
\]
and
\[
\Theta_2 = (\mu_2 +\lambda)(\alpha_2 + \beta_2 +\lambda)- \alpha_2
\left( \phi_2 - \frac{2 \phi_2 G}{c_2} + k_2 P \right).
\]

Recall that an equilibrium point ${\bf x_0}$ is said to be a {\it
saddle of type $n-p$} if the Jacobian matrix $J({\bf x_0})$ has
$n$ eigenvalues with negative real parts and $p$ eigenvalues with
positive real parts.

\begin{teo}
If $R_1 > 1$, $R_2 > 1$ and $k_1 < k_{1_{max}}$ then:
\begin{enumerate}
\item The equilibrium $\mathcal A_1$ is a saddle of type 2-2;

\item The equilibrium $\mathcal A_2$ is a saddle of type 3-1;

\item The equilibrium $\mathcal A_3$ is a saddle of type 3-1.
\end{enumerate}

\label{teo1}
\end{teo}

\noindent {\bf Proof.} From (\ref{caracteristica}) the eigenvalues
of $J(\mathcal A_1)$ are given by
\[
\begin{array}{ll}
\lambda _1 = & -\frac{1}{2}(\alpha_1 +\beta_1 +\mu_1)
+\frac{1}{2}\sqrt{(\alpha_1 +\beta_1 +\mu_1)^2 + 4\alpha_1 \phi_1
[1-\frac{1}{R_1}]},\\ \\ \lambda _2 = & -\frac{1}{2}(\alpha_1
+\beta_1 +\mu_1) -\frac{1}{2}\sqrt{(\alpha_1 +\beta_1 +\mu_1)^2 +
4\alpha_1 \phi_1 [1-\frac{1}{R_1}]},\\ \\ \lambda _3 = &
-\frac{1}{2}(\alpha_2 +\beta_2 +\mu_2 )
+\frac{1}{2}\sqrt{(\alpha_2+ \beta_2 +\mu_2 )^2 + 4 \alpha_2
\phi_2 [1-\frac{1}{R_2}]},\\ \\ \lambda _4 = &
-\frac{1}{2}(\alpha_2 +\beta_2 +\mu_2 )
-\frac{1}{2}\sqrt{(\alpha_2 +\beta_2 +\mu_2 )^2 + 4 \alpha_2 \phi_2 [1-\frac{1}{R_2}]},\end{array}\\
\]
and satisfy: $\lambda_1 > 0$, $\lambda_2 < 0$, $\lambda_3 > 0$ and
$\lambda_4 <0$. This proves the first assertion.

From (\ref{caracteristica}) the eigenvalues of $J(\mathcal A_2)$
are given by
\[
\lambda _1 = -\frac{1}{2}(\alpha_1 +\beta_1 +\mu_1 )
+\frac{1}{2}\sqrt{(\alpha_1 +\beta_1 +\mu_1 )^2 - 4 \alpha_1
\phi_1 [1-\frac{1}{R_1}]},
\]
\[
\lambda _2 = -\frac{1}{2}(\alpha_1 +\beta_1 +\mu_1 )
-\frac{1}{2}\sqrt{(\alpha_1 +\beta_1 +\mu_1 )^2 - 4 \alpha_1
\phi_1 [1-\frac{1}{R_1}]},
\]
\[
\begin{array}{cl}
\lambda _3  = & -\frac{1}{2}(\alpha_2 +\beta_2 +\mu_2 ) +
\\ \\ & \frac{1}{2}\sqrt{(\alpha_2 +\beta_2 +\mu_2 )^2 + 4 \alpha_2 \phi_2
[1-\frac{1}{R_2}] + 4 \frac{c_1 \alpha_2
\mu_1}{\alpha_1}[1-\frac{1}{R_1}]k_2},
\end{array}
\]
\[
\begin{array}{cl}
\lambda _4 = & -\frac{1}{2}(\alpha_2 +\beta_2 +\mu_2 ) -
\\ \\ & \frac{1}{2}\sqrt{(\alpha_2 +\beta_2 +\mu_2 )^2 + 4 \alpha_2 \phi_2
[1-\frac{1}{R_2}] + 4 \frac{c_1 \alpha_2
\mu_1}{\alpha_1}[1-\frac{1}{R_1}]k_2}.
\end{array}
\]
Is immediate to see that $\lambda_3 > 0$ and $\lambda_4 <0$. If
\[
\phi_1 > \frac{1}{4 \alpha_1} [(\alpha_1 + \beta_1
+\mu_1)^2+4\mu_1(\alpha_1 +\beta_1)]
\]
then $\lambda_1$ and $\lambda_2$ are complex with negative real
parts and if
\[
\phi_1 \leq \frac{1}{4 \alpha_1} [(\alpha_1 + \beta_1
+\mu_1)^2+4\mu_1(\alpha_1 +\beta_1)]
\]
then $\lambda_1 < 0$ and $\lambda_2 <0$. This proves the second
assertion.

From (\ref{caracteristica}) the eigenvalues of $J(\mathcal A_3)$
are given by
\[
\begin{array}{cl}
\lambda _1 = & -\frac{1}{2}(\alpha_1 +\beta_1 +\mu_1 +  c_2
k_1[1-\frac{1}{R_2}])+ \\ \\ & \frac{1}{2}\sqrt{(\alpha_1 +\beta_1
-\mu_1 + c_2 k_1 [1-\frac{1}{R_2}])^2 + 4\alpha_1 \phi_1},
\end{array}
\]
\[
\begin{array}{cl}
\lambda _2 = & -\frac{1}{2}(\alpha_1 +\beta_1 +\mu_1 + c_2 k_1
[1-\frac{1}{R_2}]) - \\ \\ & \frac{1}{2}\sqrt{(\alpha_1 +\beta_1
-\mu_1 + c_2 k_1 [1-\frac{1}{R_2}])^2 + 4\alpha_1 \phi_1 },
\end{array}
\]
\[
\lambda _3 = -\frac{1}{2}(\alpha_2 +\beta_2 +\mu_2 )
+\frac{1}{2}\sqrt{(\alpha_2 +\beta_2 +\mu_2 )^2 - 4\alpha_2 \phi_2
[1-\frac{1}{R_2}]},
\]
\[
\lambda _4 = -\frac{1}{2}(\alpha_2 +\beta_2 +\mu_2 )
-\frac{1}{2}\sqrt{(\alpha_2 +\beta_2 +\mu_2 )^2 - 4\alpha_2 \phi_2
[1-\frac{1}{R_2}]}.
\]
It follows that $\lambda_1 > 0$, $\lambda_2 < 0$. If
\[
\phi_2 > \frac{1}{4 \alpha_2} [(\alpha_2 + \beta_2 +\mu_2 )^2 +4
\mu_2 (\alpha_2 +\beta_2)]
\]
then $\lambda_3$ and $\lambda_4$ are complex with negative real
parts and if
\[
\phi_2 \leq \frac{1}{4 \alpha_2} [(\alpha_2 + \beta_2 +\mu_2 )^2
+4 \mu_2 (\alpha_2 +\beta_2 )]
\]
then $\lambda_3 < 0$ and $\lambda_4 <0$. This proves the last
assertion. \begin{flushright} $\blacksquare$
\end{flushright}

For the sake of completeness we state the following lemma which is
a particular case of the Theorem of Routh--Hurwitz. See
\cite{pon}, p. 62.

\begin{lema} \label{routh}
The polynomial $L(\lambda) = a_0 \lambda ^4 + a_1 \lambda ^3 + a_2
\lambda^2 + a_3 \lambda + a_4$, $a_0 > 0$, with real coefficients
has all roots with negative real parts if and only if the numbers
$a_1 , a_2 , a_3, a_4$ are positive and the inequality
\[
\Delta = a_1 \; a_2 \; a_3 - a_0 \; a_3 ^2 - a_1 ^2 \; a_4 > 0
\]
is satisfied.
\end{lema}

\begin{teo}
If $R_1 > 1$, $R_2 > 1$ and $k_1 < k_{1_{max}}$ then all the
coefficients of the characteristic polynomial of $J(\mathcal A_4)$
are positive. Therefore, if
\begin{equation} \label{critico}
\Delta = a_1 \; a_2 \; a_3 - a_3 ^2 - a_1 ^2 \; a_4
> 0,
\end{equation}
where
\[
a_1 = \alpha_1 + \beta_1 + \mu_1 + \alpha_2 + \beta_2 + \mu_2 +
k_1 G_4,
\]
\[
a_2 = \frac{\alpha_1 \phi_1}{c_1}M_4+\frac{\alpha_2 \phi_2
}{c_2}G_4 +(\alpha_1 +\beta_1 + \mu_1 +k_1 G_4)(\alpha_2 +\beta_2
+\mu_2),
\]
\[
a_3 = (\alpha_1 +\beta_1 +\mu_1 +k_1 G_4)\frac{\alpha_2 \phi_2
}{c_2} G_4 +(\alpha_2 +\beta_2 +\mu_2)\frac{\alpha_1
\phi_1}{c_1}M_4+ \alpha_2 k_1 k_2 P_4 G_4,
\]
\[
a_4 = \alpha_2 \alpha_1 \phi_1 \left( k_2 \left( 1-\frac{1}{R_1}
\right) P_4+ \frac{\phi_2}{c_1} \left(1-\frac{1}{R_2} \right) M_4
\right),
\]
then the differential equations (\ref{sistemafinal}) have an
asymptotically stable equilibrium point at $\mathcal A_4$. If
\[
\Delta < 0
\]
then $\mathcal A_4$ is unstable.
\label{teo2}
\end{teo}

\noindent {\bf Proof.} From (\ref{caracteristica})
the characteristic polynomial of $J(\mathcal A_4)$ is given by
\[ [\lambda ^2+(\alpha_1 +\beta_1 +\mu_1 +k_1 G_4)\lambda +\frac{\alpha_1 \phi_1}{c_1}M_4]
[\lambda ^2+(\alpha_2 +\beta_2 +\mu_2)\lambda +\frac{\alpha_2
\phi_2}{c_2} G_4]\]
\[ + \alpha_2 \mu_1 k_1 k_2 P_4 G_4 + \alpha_2 k_1 k_2 P_4 G_4 \lambda,\]
which can be written as
\[ \begin{array}{l}
\lambda ^4 +\lambda ^3 [\alpha_1 +\beta_1 +\mu_1 +\alpha_2
+\beta_2+ \mu_2 +k_1 G_4] \\ \\ +\lambda ^2 [\frac{\alpha_1
\phi_1}{c_1} M_4+\frac{\alpha_2 \phi_2 }{c_2} G_4 +(\alpha_1
+\beta_1 +\mu_1 +k_1 G_4)(\alpha_2 +\beta_2 +\mu_2)]\\ \\+\lambda
[(\alpha_1 +\beta_1 +\mu_1 +k_1 G_4)\frac{\alpha_2 \phi_2}{c_2}
G_4 +(\alpha_2 +\beta_2 +\mu_2)\frac{\alpha_1 \phi_1}{c_1}M_4+
\alpha_2 k_1 k_2 P_4 G_4]\\ \\+\frac{\alpha_1 \phi_1}{c_1}M_4
\frac{\alpha_2 \phi_2}{c_2} G_4 + \alpha_2 \mu_1 k_1 k_2 P_4 G_4.
\end{array} \]
Now it is simple to see that the coefficients of the
characteristic polynomial are given by $a_1, a_2, a_3, a_4$ above.
From the hypotheses these coefficients are positive. The stability
at $\mathcal A_4$ follows from Lemma \ref{routh}.
\begin{flushright}
$\blacksquare$
\end{flushright}

The following corollary is immediate from the fact that $a_i >0$.

\begin{cor}
The Jacobian matrix $J(\mathcal A_4)$ has a pair of complex
eigenvalues with zero real part if and only if
\begin{equation}
a_3 ^2 - a_1 \; a_2 \; a_3 + a_1 ^2 \; a_4 =0, \label{condhopf}
\end{equation}
where $a_i$ are defined in Theorem \ref{teo2}.

\label{hopf}
\end{cor}

In next section we study the stability of $\mathcal A_4$ under the
condition (\ref{condhopf}), complementary to the range of validity
of Theorem \ref{teo2}.

\section{Hopf Bifurcation Analysis}\label{hopfbif}

\subsection{Generalities on Hopf Bifurcations}\label{generalities}

The study outlined below is based on the approach found in the book
of Kuznetsov \cite{kuznet}, pp 175-178.

Consider the differential equations
\begin{equation}
{\bf x}' = f ({\bf x}, {\bf \mu}), \label{diffequat}
\end{equation}
\noindent where ${\bf x} \in \R^4$ and ${\bf \mu} \in \R^m$ is a
vector of control parameters. Suppose (\ref{diffequat}) has an
equilibrium point ${\bf x} = {\bf x_0}$ at ${\bf \mu} = {\bf
\mu_0}$ and represent
\begin{equation}
F({\bf x}) = f ({\bf x}, {\bf \mu_0}) \label{Fhomo}
\end{equation}
as
\begin{equation}
F({\bf x}) = A{\bf x} + \frac{1}{2} \: B({\bf x},{\bf x}) +
\frac{1}{6} \: C({\bf x}, {\bf x}, {\bf x}) + O(|| {\bf x}
||^4){\nonumber}, \label{taylorexp}
\end{equation}
\noindent where $A = f_{\bf x}(0,{\bf \mu_0})$ and
\begin{equation}
B_i ({\bf x},{\bf y}) = \sum_{j,k=1}^4 \frac{\partial ^2
F_i(\xi)}{\partial \xi_j \: \partial \xi_k} \bigg|_{\xi=0} x_j \;
y_k, \label{Bap}
\end{equation}
\begin{equation}
C_i ({\bf x},{\bf y},{\bf z}) = \sum_{j,k,l=1}^4 \frac{\partial ^3
F_i(\xi)}{\partial \xi_j \: \partial \xi_k \: \partial \xi_l}
\bigg|_{\xi=0} x_j \; y_k \: z_l, \label{Cap}
\end{equation}
\noindent for $i = 1, 2, 3, 4$. Here the variable ${\bf x}-{\bf
x_0}$ is also denoted by ${\bf x}$.

Suppose $({\bf x_0}, {\bf \mu_0})$ is an equilibrium point of
(\ref{diffequat}) where the Jacobian matrix $A$ has a pair of
purely imaginary eigenvalues $\lambda_{3,4} = \pm i \omega_0$,
$\omega_0 > 0$, and no other critical (i.e., on the imaginary
axis) eigenvalues.

Let $p, q \in \C ^4$ be vectors such that
\begin{equation}
A q = i \omega_0 \: q,\:\: A^{\top} p = -i \omega_0 \: p, \:\:
\langle p,q \rangle = \sum_{i=1}^4 \bar{p}_i \: q_i \:\: = 1.
\label{normalization}
\end{equation}

The two dimensional center manifold can be parameterized by $w \in
\R^2 = \C$,  by means of ${\bf x} = H (w,\bar w )$, which is
written as
\[
H(w,{\bar w}) = w q + {\bar w}{\bar q} + \sum_{2 \leq j+k \leq 3}
\frac{1}{j!k!} \: h_{jk}w^j{\bar w}^k + O(|w|^4),
\]
with $h_{jk} \in \C ^4$, $h_{jk}={\bar h}_{kj}$.

Substituting these expressions into (\ref{diffequat}) and
(\ref{taylorexp}) one has
\begin{equation}
H_w (w,\bar w )w' + H_{\bar w}  (w,\bar w ){\bar w}'  = F(H(w
,\bar w )). \label{homologicalp}
\end{equation}

The complex vectors $h_{ij}$ are to be determined so that equation
(\ref{homologicalp}) writes as follows
\[
w'= i \omega_0 w + \frac{1}{2} \: G_{21} w |w|^2 + O(|w|^4),
\]
with  $G_{21} \in \C $.

Solving the linear system obtained by expanding
(\ref{homologicalp}), the coefficients of the quadratic terms of
(\ref{Fhomo}) lead to
\begin{equation}
h_{11}=-A^{-1}B(q,{\bar q}) \label{h11},
\end{equation}
\begin{equation}
h_{20}=(2i\omega_0 I_4 - A)^{-1}B(q,q),\label{h20}
\end{equation}
where $I_4$ is the unit $4 \times 4$ matrix.

The coefficients of the cubic terms are also uniquely calculated,
except for the term $w^2 {\bar w}$, whose coefficient satisfies a
singular system for $h_{21}$
\begin{equation}
(i \omega_0 I_4 -A)h_{21}=C(q,q,{\bar q})+B({\bar q},h_{20}) + 2
B(q,h_{11})-G_{21}q, \label{h21m}
\end{equation}
which has a solution if and only if
\[
\langle p, C(q,q,\bar q) + B(\bar q, h_{20}) + 2 B(q,h_{11})
-G_{21} q \rangle = 0.
\]
Therefore
\begin{equation} \label{G21}
G_{21}= \langle p, C(q,q,\bar q) + B(\bar q, (2i \omega_0
I_4-A)^{-1} B(q,q)) - 2 B(q,A^{-1} B(q,\bar q)) \rangle,
\end{equation}
and the {\it first Lyapunov coefficient} $l_1$ -- which decides by
the analysis of third order terms at the equilibrium its
stability, if negative, or instability, if positive -- is defined
by
\begin{equation}
l_1 =  \frac{1}{2 \; \omega_0} \: {\rm Re} \; G_{21}.
\label{defcoef}
\end{equation}

A {\it Hopf point} $({\bf x_0}, {\bf \mu_0})$ is an equilibrium
point of (\ref{diffequat}) where the Jacobian matrix $A$ has a
pair of purely imaginary eigenvalues $\lambda_{3,4} = \pm i
\omega_0$, $\omega_0 > 0$, and no other critical eigenvalues. At a
Hopf point, a two dimensional center manifold is well-defined,
which is invariant under the flow generated by (\ref{diffequat})
and can be smoothly continued to nearby parameter values.

A Hopf point is called {\it transversal} if the curves of complex
eigenvalues cross the imaginary axis with non-zero derivative.

In a neighborhood of a transversal Hopf point with $l_1 \neq 0$
the dynamic behavior of the system (\ref{diffequat}), reduced to
the family of parameter-dependent continuations of the center
manifold, is orbitally topologically equivalent to the complex
normal form
\begin{equation}\label{nf}
w' = (\gamma + i \omega) w + l_1 w |w|^2 ,
\end{equation}
$w \in \C $, $\gamma$, $\omega$ and $l_1$ are smooth continuations
of $0$, $\omega_0$ and the first Lyapunov coefficient at the Hopf
point \cite{kuznet}, respectively. When $l_1 < 0$ ($l_1
> 0$) a family of stable (unstable) periodic orbits can be found
on this family of center manifolds, shrinking to the equilibrium
point at the Hopf point.

\subsection{Hopf Bifurcation in the Biological Model} \label{ss:hopfbif}

In this subsection we analyze the stability at $\mathcal A_4$
given by (\ref{ponto4}) under the condition (\ref{condhopf}). From
(\ref{diffequat}) write the Taylor expansion (\ref{taylorexp}) of
$f({\bf x})$. Thus
\begin{equation}\label{A}
A = \left(\begin{array}{cccc} -(\alpha_1 + \beta_1) -k_1 G_4 & \phi_1 (1-\frac{2 M_4}{c_1}) & 0 & -k_1 P_4 \\
\alpha_1 & - \mu_1 & 0 & 0
\\ k_2 G_4 & 0 & - (\alpha_2 + \beta_2 ) & \phi_2 (1-\frac{2 G_4}{c_2})+k_2 P_4 \\
0 & 0 & \alpha_2 & - \mu_2 \end{array} \right)
\end{equation}
and, with the notation in (\ref{taylorexp}) to (\ref{Cap}), one
has
\begin{equation}
F({\bf x})\, - \, A{\bf x} = \left( - \frac{\phi_1 M^2}{c_1} - k_1
P G, 0, - \frac{\phi_2 G^2}{c_2} +  k_2 P G ,0  \right).
\label{partenaolinear}
\end{equation}

From (\ref{taylorexp}), (\ref{Bap}), (\ref{Cap}) and
(\ref{partenaolinear}) one has
\begin{equation}
B({\bf x},{\bf y}) = \left( B_1({\bf x},{\bf y}),0, B_3({\bf
x},{\bf y}), 0 \right), \label{B1}
\end{equation}
where
\[
B_1({\bf x},{\bf y})= -\frac{2\phi_1}{c_1} \: x_2 \; y_2 - k_1
(x_1 \; y_4 + x_4 \; y_1),
\]
\[
B_3({\bf x},{\bf y}) = -\frac{2 \phi_2}{c_2} \: x_4 \; y_4 + k_2
(x_1 \; y_4 + x_4 \; y_1),
\]
and
\begin{equation}
C({\bf x},{\bf y}, {\bf z}) \equiv 0. \label{C1}
\end{equation}

To pursue the analysis  consider the following table of specific
parameters
\begin{equation}
\label{basedados}
\begin{array}{|lllll|}
\hline \alpha_1 = 0.7 & \beta_1 =0.003 &\mu_1 = 0.6 & \phi_1= 2.3
&
    c_1 = 400000  \\ \alpha_2 =0.3 & \beta_2 =0.0015 & \mu_2 =0.4  & \phi_2 = 4 &  c_2
    = 100 \\ \hline
    \end{array}
\end{equation}
taken from \cite{dan} and \cite{son}, where their biological
feasibility in Brazilian fields is discussed.

With the above parameter values the differential equations
(\ref{sistemafinal}) are in fact a two parameter system of
differential equations where the parameters are $k_1$ and $k_2$
and can be written equivalently as
\begin{equation}
{\bf x}' = f({\bf x}, k_1, k_2), \label{campo}
\end{equation}
with $f({\bf x}, k_1, k_2)$ defined by the right-hand sides of
(\ref{sistemafinal}).

\begin{figure}[!h]
\centerline{
\includegraphics[width=8cm]{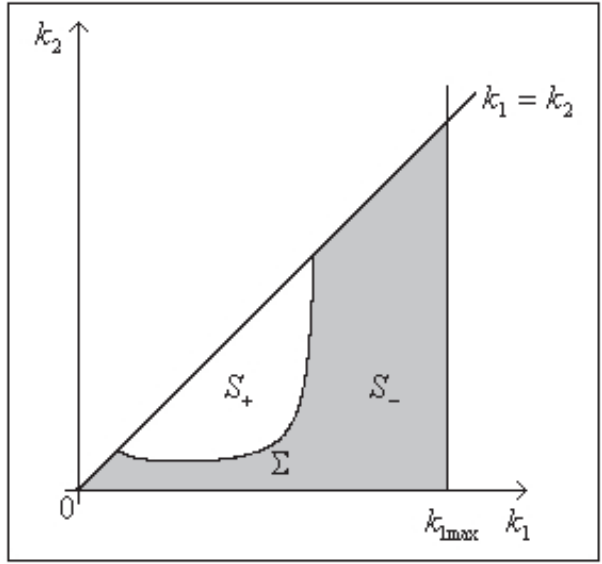}}
\caption{{\small Set of admissible parameters $\mathcal S$ and
Hopf curve $\Sigma$}.}

\label{admpara}
\end{figure}

With the parameter values of table (\ref{basedados}), the
equilibrium point $\mathcal A_4$ (\ref{ponto4}) has the following
coordinates
\[
P_4 = \frac{800000 \; (1.425 - 64.764 \; k_1)}{4.508 + 1.444 \cdot
10^7 \; k_1 \; k_2}, \: M_4 = \frac{933333.333 \; (1.425 - 64.764
\; k_1)}{4.508 + 1.444 \cdot 10^7 \; k_1 \; k_2},
\]
\[
L_4 = \frac{444.444 \; (1.216 + 85550.400 \; k_2)}{4.508 + 1.444
\cdot 10^7 \; k_1 \; k_2}, \: G_4 = \frac{333.333 \; (1.216 +
85550.400 \; k_2)}{4.508 + 1.444 \cdot 10^7 \; k_1 \; k_2},
\]
while $R_1$, $R_2$ and $k_{1_{max}}$, given by (\ref{R1R2}) and
(\ref{k1max}), have the form
\begin{equation}
R_1 = 3.81697, \: \: R_2 = 9.95025, \: \: k_{1_{max}} = 0.0220159.
\label{valoresR1R2k1max}
\end{equation}
From the above equation and the Remark \ref{parametros}, the set
of admissible parameters is given by (see Fig \ref{admpara})
\begin{equation}
\mathcal S = \{(k_1,k_2) | \: 0 < k_1 < k_{1_{max}} = 0.0220159
\:\: {\mbox {and}} \: \: 0 < k_2 \leq k_1 \}.
 \label{admissible}
\end{equation}

In this set $\mathcal S$  the curve $\Sigma = \Delta^{-1} (0)$ is
well-defined (see (\ref{condhopf})), where $\Delta$ is given by
\begin{eqnarray*}
1699.422 - 233762.372 k_1 - 6.860 \cdot 10^7 k_2 - 1.114 \cdot
10^7 k_1^2 - 2.175 \cdot 10^{12} k_2^2 - \\ 4.994 \cdot 10^{10}
k_1 k_2 - 2.147 \cdot 10^8 k_1^3 - 4.079 \cdot 10^{12} k_1^2 k_2 -
1.529 \cdot 10^{15} k_1 k_2^2 - \\ 7.809 \cdot 10^{13} k_1^3 k_2 -
4.319 \cdot 10^{17} k_1^2 k_2^2 - 1.540 \cdot 10^{19} k_1 k_2^3 -
4.755 \cdot 10^{14} k_1^4 k_2 - \\ 1.752 \cdot 10^{19} k_1^3 k_2^2
- 6.741 \cdot 10^{21} k_1^2 k_2^3 + 2.940 \cdot 10^{18} k_1^4
k_2^2 - 1.703 \cdot 10^{24} k_1^3 k_2^3 + \\ 1.634 \cdot 10^{26}
k_1^2 k_2^4 + 6.618 \cdot 10^{24} k_1^4 k_2^3 - 4.437 \cdot
10^{28} k_1^3 k_2^4 + 1.643 \cdot 10^{29} k_1^4 k_2^4,
\end{eqnarray*}
representing the parameters where $J(\mathcal A_4)$ has a pair of
purely imaginary eigenvalues $\lambda_{3,4} = \pm i \omega_0$ with
\begin{eqnarray}\label{omega0}
\omega_0 = 1.2909 \Big [ 0.6909 + \frac{0.0071 (-2.6479 \cdot
10^{-6} + k_2)(1.4219 \cdot 10^{-5} + k_2)}{k_2 (3.1305 \cdot
10^{-7} + k_1
k_2)} + \nonumber \\
 \frac{8.5745 \cdot 10^{-7}}{k_2} + \Big( \frac{1}{(6.2611 \cdot
10^{-7} + k_1 k_2)k_1 k_2} \Big( (1.0783 \cdot 10^{-9} + \\
5.0825 \cdot 10^{-5} k_2) k_2 - 3.1286 \cdot 10^{-4} (1.3649 \cdot
10^{-3} k_2) (6.0127 \cdot 10^{-6} + k_2) + \nonumber \\
0.9391 \; k_1 ^2 (9.5980 \cdot 10^{-8} + k_2)(8.1563 \cdot 10^{-6}
+ k_2) \Big) \Big )^{1/2} \Big ]^{1/2}. \nonumber
\end{eqnarray}

Thus one has (see Fig. \ref{admpara})
\[
\mathcal S = \mathcal S_{+}  \cup \Sigma \cup \mathcal S_{-}.
\]
For parameter values in the region $\mathcal S_{+}$ the
equilibrium $\mathcal A_4$ is unstable since the Jacobian matrix
$J(\mathcal A_4)$ has two complex eigenvalues with positive real
parts and two other real negative eigenvalues. For parameter
values in the region $\mathcal S_{-}$ the equilibrium $\mathcal
A_4$ is asymptotically stable since $J(\mathcal A_4)$ has two
complex eigenvalues with negative real parts and two other real
negative eigenvalues. The curve $\Sigma$ is the curve of
admissible parameters where the equilibrium $\mathcal A_4$ is a
Hopf point.

\begin{teo}
Consider the differential equations (\ref{sistemafinal}) with the
parameters given  in the table (\ref{basedados}). If $(k_1,k_2)
\in \Sigma$ then the two parameter family of differential
equations (\ref{sistemafinal}) has a transversal Hopf point at
$\mathcal A_4$. This Hopf point at $\mathcal A_4$ is unstable and
for each $(k_1,k_2) \in \mathcal S_{-}$, but close to $\Sigma$,
there exists an unstable periodic orbit near the asymptotically
stable equilibrium point $\mathcal A_4$. See Fig \ref{admpara}.

\label{teoremahopf}
\end{teo}

\noindent {\bf Computer Assisted Proof.} The proof follows the
steps outlined in  Subsection \ref{generalities}. However, all the
expressions in the proof are too long to be put in print. For this
reason, in the site \cite{mello} have been posted the main steps
of the long calculations involved in the proof. This has been done
in the form of a {\it notebook} for MATHEMATICA 5 \cite{math}. A
sufficient condition for being a Hopf point is that the first
Lyapunov coefficient $l_1 \neq 0$. In fact, it can be shown
numerically that $l_1 (k_1, k_2) > 0$ for all values $(k_1, k_2)
\in \Sigma$. A particular case and other related calculations are
considered below for the sake of illustration.

Take the particular point $Q =(k_1 = 0.00331, k_2 = 0.00100) \in
\Sigma $ with five decimal round-off coordinates \cite{dan}. For
these values of the parameters
\[
\mathcal A_4 = (18543.57758, 21634.17385, 738.0525862,
553.5394397).
\]
The Jacobian matrix $J(\mathcal A_4)$ has eigenvalues
\[
\lambda_1 = - 3.61058, \: \lambda_2 = - 0.22912, \: \lambda_{3,4}
= \pm 2.84670 i,
\]
and thus
\begin{equation}
\omega_0 = 2.84670.
\label{omegazeronum}
\end{equation}
From (\ref{normalization}) the eigenvectors $q$ and $p$ have the
form
\[
q= \left( \begin{array}{c} 820.5542609+1080.774610i \\
    295.1756045-139.5588184i \\ 862.8021803+130.4940530i \\
    26.01486634-87.27100717i \end{array} \right),
\]
\[
p = \left( \begin{array}{c} 0.00003314748646+0.00006274424412i \\
    -0.00003846764141+0.00003199241887i \\ 0.0005233172006+0.00007678211168i \\
     0.001254520214-0.004888597529i \end{array} \right).
\]
One has
\[
B(q,q)=\left( \begin{array}{c}-767.7418261+289.3499796i \\ 0 \\
  786.4902302+276.2661945i \\ 0 \end{array} \right)
\]
and
\[
B(q,\bar{q})=\left( \begin{array}{c} 482.6477605 \\ 0 \\
    -809.3875158+0.6\cdot10^{-7}i \\ 0 \end{array}\right).
\]
From (\ref{h11}) and (\ref{h20}) the complex vectors $h_{11}$ and
$h_{20}$ have the form
\[
- h_{11}= \left( \begin{array}{c} -1622.977370+0.9904140546 \cdot 10^{-7}i \\
   -1893.473598+0.1155483063 \cdot 10^{-6}i \\
 -5.359116331-0.3117318364 \cdot 10^{-9}i \\
  -4.019337253-0.2337988771 \cdot 10^{-9}i \end{array} \right),
\]
\[
h_{20}= \left( \begin{array}{c} 71.87520338+68.12253398i \\
    9.204672581-7.866965075i \\ 83.57142169-174.4554220i \\
      -8.839482676-5.024621730i \end{array} \right).
\]
From (\ref{G21}) the complex number $G_{21}$ is given by
\begin{equation}
G_{21}= 0.057297 - 0.027485 i, \label{G21numerico}
\end{equation}
and from (\ref{defcoef}), (\ref{omegazeronum}) and
(\ref{G21numerico}) one has the first Lyapunov coefficient at $Q$
\begin{equation}
l_1 (Q) = 0.00353522 > 0. \label{L1Q}
\end{equation}
The above calculations have also been checked with 10 decimals
round-off precision performed using the software MATHEMATICA 5
\cite{math}. See \cite{mello}.

Some other values of pairs $(k_1,k_2) \in \Sigma$, the values of
the complex eigenvalues of $J(\mathcal A_4)$ as well as the
corresponding values of $l_1 (k_1,k_2)$ are listed the table
below. The calculations leading to these values can be found in
\cite{mello}.

\vspace{1cm}
\begin{tabular}{|c|c|c|c|}    \hline \hline
$k_1$       & $k_2$     & complex eigenvalues of $J(\mathcal A_4)$  &$l_1 (k_1,k_2)$ \\
\hline 0.0004813    & 0.0004812   &   $\pm 4.76456 i$     &
$4.69457 \cdot 10^{-8}$ \\ \hline
 0.0007954    & 0.0003535   & $\pm 3.98051 i$       & $1.21597 \cdot
 10^{-7}$ \\ \hline
 0.0011096  & 0.0003086   & $\pm 3.59051 i $      & $2.17614 \cdot 10^{-7}$ \\ \hline
 0.0014238  & 0.0002950   & $\pm 3.35518 i $      & $3.29018 \cdot 10^{-7}$ \\ \hline
 0.0017379  & 0.0003001   & $\pm 3.19780 i $      & $4.52683 \cdot 10^{-7}$ \\ \hline
 0.0020521  & 0.0003220   & $\pm 3.08560 i $      & $5.86835 \cdot 10^{-7}$ \\ \hline
 0.0023663  & 0.0003649   & $\pm 3.00207 i $      & $7.30346 \cdot 10^{-7}$ \\ \hline
 0.0026804  & 0.0004427   & $\pm 2.93795 i $      & $8.82421 \cdot 10^{-7}$ \\ \hline
 0.0029946  & 0.0005957   & $\pm 2.88762 i $      & $1.04245 \cdot 10^{-6}$ \\ \hline
 0.0033088  & 0.0009855   & $\pm 2.84745 i $      & $1.20994 \cdot 10^{-6}$ \\ \hline
 0.0036230  & 0.0035924   & $\pm 2.81501 i $      & $1.38449 \cdot 10^{-6}$ \\ \hline\hline
\end{tabular}

\begin{flushright} $\blacksquare$ \end{flushright}

\begin{remark}
The value of the first Lyapunov coefficient $l_1$ does depend on the
normalization of the eigenvectors $q$ and $p$, while its sign is
invariant under scaling of $q$ and $p$ obeying the relative
normalization. See \cite{kuznet}, p. 99. The values $l_1(Q)$ in
(\ref{L1Q}) and $l_1(k_1,k_2)$ in the above Table are obtained with
two different choices of the eigenvectors $q$ and $p$, see
\cite{mello}. This explains the difference in the order of magnitude
of the numbers involved.

\end{remark}

\begin{figure}[!h]
\centerline{
\includegraphics[width=9cm]{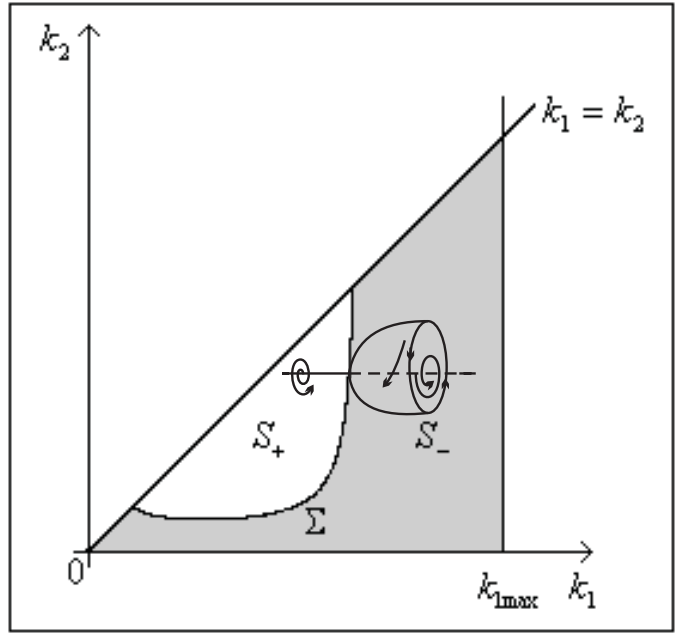}}
\caption{{\small Bifurcation diagram for a typical point on the
curve $\Sigma$}.}

\label{fig3}
\end{figure}

As a consequence of Theorem 3.1 there are no Hopf points of
codimension 2 on $\Sigma$ since the sign of the first Lyapunov
coefficient does not change.
In Fig. \ref{fig3} is illustrated the bifurcation diagram for a
typical point on the curve $\Sigma$.

Assuming the same values in the table (\ref{basedados}) in next
theorem we study the behavior of the Hopf curve $\Sigma$ in the set
of admissible parameters $\mathcal S$ (see equation
(\ref{admissible})) as the parameter $c_2$ increases. In fact, the
carrying capacity, representing several other factors, has a
determinant role on the populations under study.

\begin{figure}[!h]
\centerline{
\includegraphics[width=7cm]{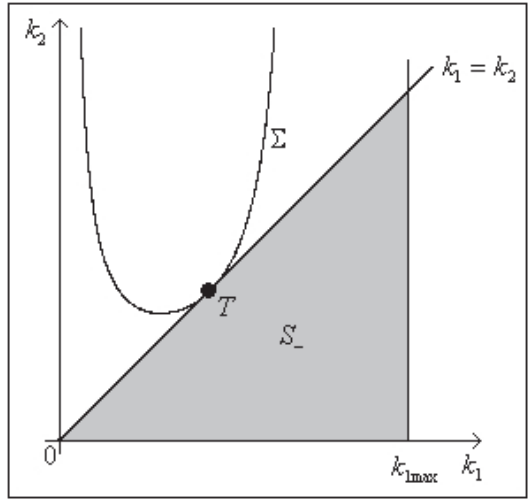}}
\caption{{\small Curve $\Sigma$ intersects $\mathcal S$ at one
point $T$}.}

\label{curvasigma}
\end{figure}

\begin{teo}
The one parameter family of curves $\Sigma_{c_2}=\Delta_{c_2}^{-1}
(0)$ has only one point of tangency $T$ with the line $k_1=k_2$
for $c_2 = 650.41463$. For values $c_2 > 650.41463$ the curve
$\Sigma_{c_2}$ does not intersect the set $\mathcal S$. Therefore
for values $c_2 > 650.41463$ the set $\mathcal S_{+}$ is empty,
$\mathcal S = \mathcal S_{-}$ and the equilibrium $\mathcal A_4$
is asymptotically stable for all values $(k_1, k_2) \in \mathcal
S$. See Fig. \ref{curvasigma}.

\label{moves}
\end{teo}

\noindent {\bf Proof.} The surface of Hopf points, or equivalently
the one parameter family of Hopf curves, where $J(\mathcal A_4)$
has a pair of purely imaginary eigenvalues is defined by
$\Sigma_{c_2} = \{ \Delta(k_1,k_2,c_2) = 0 \}$ where
$\Delta(k_1,k_2,c_2)$ (see (\ref{condhopf})) is given by
\begin{eqnarray*}
1699.422 - 2337.623 c_2 k_1 - 6.860 \cdot 10^7 k_2 - 1114.941
c_2^2 k_1^2 - 2.175 \cdot 10^{12} k_2^2 - \\ 4.994 \cdot 10^{8}
c_2 k_1 k_2 - 214.747 c_2^3 k_1^3 - 4.079 \cdot 10^{8} c_2^2 k_1^2
k_2 - 1.529 \cdot 10^{13} c_2 k_1 k_2^2 - \\ 7.809 \cdot 10^{7}
c_2^3 k_1^3 k_2 - 4.319 \cdot 10^{13} c_2^2 k_1^2 k_2^2 - 1.540
\cdot 10^{17} c_2 k_1 k_2^3 - \\ 4.755 \cdot 10^{6} c_2^4 k_1^4
k_2 - 1.752 \cdot 10^{13} c_2^3 k_1^3 k_2^2 - 6.741 \cdot 10^{17}
c_2^2 k_1^2 k_2^3 + \\ 2.940 \cdot 10^{10} c_2^4 k_1^4 k_2^2 -
1.703 \cdot 10^{18} c_2^3 k_1^3 k_2^3 + 1.634 \cdot 10^{22} c_2^2
k_1^2 k_2^4 + \\ 6.618 \cdot 10^{16} c_2^4 k_1^4 k_2^3 - 4.437
\cdot 10^{22} c_2^3 k_1^3 k_2^4 + 1.643 \cdot 10^{21} c_2^4 k_1^4
k_2^4.
\end{eqnarray*}
The intersection of the surface $\Sigma_{c_2}$ with the plane $k_1
= k_2$ determines the curve $\mathcal C$, given implicitly by
\begin{eqnarray*}
N(k_1,c_2) = 1699.422 - ( 2337.623 c_2 + 6.860 \cdot 10^7) k_1 - (
1114.941 c_2^2 + \\ 2.175 \cdot 10^{12} + 4.994 \cdot 10^{8} c_2 )
k_1^2 - ( 214.747 c_2^3 + 4.079 \cdot 10^{8} c_2^2 + \\ 1.529
\cdot 10^{13} c_2 ) k_1^3 - (7.809 \cdot 10^{7} c_2^3 + 4.319
\cdot 10^{13} c_2^2 + 1.540 \cdot 10^{17} c_2 ) k_1^4 - \\ ( 4.755
\cdot 10^{6} c_2^4 + 1.752 \cdot 10^{13} c_2^3 + 6.741 \cdot
10^{17} c_2^2 ) k_1^5 + ( 2.940 \cdot 10^{10} c_2^4 - \\ 1.703
\cdot 10^{18} c_2^3 + 1.634 \cdot 10^{22} c_2^2 ) k_1^6 + ( 6.618
\cdot 10^{16} c_2^4 - 4.437 \cdot 10^{22} c_2^3 ) k_1^7 + \\ 1.643
\cdot 10^{21} c_2^4 k_1^8 = 0.
\end{eqnarray*}

Differentiating implicitly the above expression with respect to
$k_1$ one has
\[
\frac{d c_2}{d k_1} = - \frac{ \frac {\partial N}{\partial
k_1}}{\frac{\partial N}{\partial c_2}} = 0, \: \: \frac{d^2 c_2}{d
k_1 ^2} < 0,
\]
at $k_1 = 0.00035$ and $c_2 = 650.41463$. Therefore the curve
$\mathcal C$ is a graph near the point $(k_1 = 0.00035, c_2 =
650.41463)$ and has a local maximum point at $k_1 = 0.00035 $. It
can be shown \cite{mello} that this maximum is global since ${d
c_2}/{d k_1}$ has no other zero. It is easy to verify through a
calculation that the point $T = (k_1,k_2)=(0.00035,0.00035)$
belongs to $\Sigma_{c_2}$ for $c_2 = 650.41463$. Now the gradient
of $\Delta_{c_2}$ at $T$ for $c_2 = 650.41463$ is given by
\[
(-1.73746 \cdot 10^{10}, 1.73746 \cdot 10^{10}),
\]
which is parallel to the vector $(-1,1)$, the normal to the line
$k_1=k_2$.
\begin{flushright}
$\blacksquare$
\end{flushright}

\begin{remark}
Since $k_{1_{max}}$ does depend
on the parameter $c_2$, according to Eq. (\ref{k1max}), so does the
admissible region $\mathcal S = \mathcal S_{c_2}$. For $c_2 =
650.41463$, a calculation gives $k_{1_{max}}= 0.00338491$. This is
compatible with position of $T$ at $k_1 = k_2 = 0.00035$, as
illustrated in Figure \ref{curvasigma}.

\end{remark}

\section{Concluding Comments}\label{conclusion}

In this paper we studied the system (\ref{sistemafinal}) of
interest as a mathematical model for biological control, proposed
by Yang and Ternes \cite{hyu, son, hyuson}  and studied also by
Santos \cite{dan}. Valuable field data are provided in \cite{son},
valid for the {\it citrus leafminer} and its native and imported
enemies in the region of Limeira (S\~ao Paulo, Brazil). An
extensive, enlightening discussion of the economic and
agricultural interest of the problem, other pertinent differential
equations models as well as extensive bibliography, are also
presented there.

Under conditions made explicit in Remark \ref{remark1} we
determine the unique equilibrium point ($\mathcal A_4 $) with
positive coordinates and establish necessary and sufficient
conditions for its (Lyapunov) stability (Theorem \ref{teo2}). It
can be seen however that this condition $\Delta > 0$, when
expressed in terms of the parameters is a rational function whose
denominator does not vanish and its numerator is a polynomial of
too many terms to be put in print, but still amenable to numerical
calculations. For this reason the treatment of the stability of
($\mathcal A_4 $) in Subsection \ref{ss:hopfbif} is computer
assisted. The conclusion of this study, made precise in Theorem
\ref{teoremahopf}, is the existence of periodic orbits obtained by
Hopf bifurcation, on the side (of $\Delta = 0$) where $\mathcal
A_4 $ is an attractor.

The study of the general analytic and geometric properties of the
boundary of the stability region, given by the Hopf variety
$\Delta = 0$, so as to include parameter values of biological
interest as proposed here as well as others appearing in the work
of Ternes \cite{son}, remain at the present moment as a
mathematical challenge. Theorem \ref{moves} gives only a thin
slice of the geometry.

The reports in \cite{fao} and \cite{flo}, among many others, show
that the interest for the combat of the {\it citrus leafminer}
extends to most regions where citrus trees grow.

The  Mathematica notebooks \cite{mello}, with the table
(\ref{basedados}), used in the computer assisted arguments for the
proofs of Theorems \ref{teoremahopf} and \ref{moves}, can be
adapted to tables with data pertinent to other geographic and
climatic regions and involving different host--parasitoid
interactions.

In \cite{hyuson} Ternes and Yang discuss, with pertinent
documentation, the introduction of a foreign parasitoid, {\it
Ageniaspis citricola} to add the native {\it Galeopsomya Fausta}
in the combat with the {\it leafminer,  Phyllocnistis citrella}.
They propose a model with eight differential equations for the
three species and their immature stages. In \cite{hyuson} and
\cite{son} are given starting steps for an analysis of the
stability of the equilibria in this extended eight -- dimensional
system. Based in a numerical study of a complex equilibrium point
they recommend that the biological control of the {\it leafminer}
be implemented with both the native and foreign parasitoids.
Meanwhile, the Hopf bifurcation analytic and computer algebra
study of the complex equilibria of the eight equations, with the
methods used in the present paper, seems unsurmountable at the
present moment, due to the large number of parameters involved.

\vspace{0.3cm} \noindent {\bf Acknowledgement}: The first and second
authors developed this work under the projects CNPq Grants
473824/04-3 and 473747/2006-5. The first author is fellow of CNPq.
The fourth author is supported by CAPES. This work was finished
while the second author visited Universitat Aut\`onoma de Barcelona,
supported by CNPq grant 210056/2006-1.

\end{document}